\documentclass[preprint,12pt]{elsarticle}

\usepackage{amssymb}
\usepackage{mathrsfs}
\usepackage[boxed,commentsnumbered]{algorithm2e}

\usepackage{amsmath}
\usepackage[T1]{fontenc}
\usepackage[latin9]{inputenc}
\usepackage{graphicx}
\usepackage{amssymb}
\usepackage{lscape}

\newtheorem{theorem}{Theorem}[section]
\newtheorem{lemma}[theorem]{Lemma}
\newtheorem{proposition}[theorem]{Proposition}
\newtheorem{corollary}[theorem]{Corollary}
\newtheorem{example}[theorem]{Example}
\newtheorem{remark}[theorem]{Remark}
\newtheorem{notation}[theorem]{Notation}
\newtheorem{definition}[theorem]{Definition}
\newenvironment{proof}[1][Proof]{\begin{trivlist} \item[\hskip \labelsep {\bfseries #1}]}{\end{trivlist}}

\journal{...}

\begin{document}

\begin{frontmatter}



\title{Reduced Collatz Dynamics is Periodical and the Period Equals 2 to the Power of the Count of x/2}


\author{Wei Ren} \ead{weirencs@cug.edu.cn}
\address{School of Computer Science \\ China University of
Geosciences, Wuhan, China}


\begin{abstract}

In this paper, we prove that reduced dynamics on Collatz conjecture
is periodical, and its period equals 2 to the power of the count of
x/2 computation in the reduced dynamics. More specifically, if there
exists reduced dynamics of x (that is, start from an integer x and
the computation will go to an integer less than x), then there must
also exist reduced dynamics of x+P (that is, if starting from an
integer x+P, then computation will go to an integer less than x+P),
where P equals 2 to the power of L, and L is the total count of x/2
computations (i.e., computational times) in the reduced dynamics of
x (note that, equivalently, L is also the length of the reduced
dynamics of x). Therefore, the power (or output) of this period
property, which is discovered and proved in this paper, is - the
study of the existence of reduced dynamics of x will result in the
existence of reduced dynamics of x+P (and iteratively x+n*P, n is a
positive integer). Hence, only partition of integers needs to be
verified for the existence of their reduced dynamics. Finally, if
any starting integer x can be verified for the existence of its
reduced dynamics, then Collatz Conjecture will be True (due to our
proposed Reduced Collatz Conjecture).

\end{abstract}


\end{frontmatter}





The Collatz conjecture can be stated simply as follows: Take any
positive integer $x$. If $x$ is even, divide it by $2$ to get $x/2$.
If $x$ is odd, multiply it by $3$ and add $1$ to get $3x + 1$.
Repeat the process again and again. The Collatz conjecture is that
no matter what the integer (i.e., $x$) is taken, the process will
always eventually reach 1.

The current known integers that have been verified are about 60 bits
by T.O. Silva using normal personal computers
\cite{UpperboundRecord1,UpperboundRecord2}. They verified all
integers that are less than 60 bits. Wei Ren et al. \cite{weiuic}
verified $2^{100000}-1$ can return to 1 after 481603 times of $3x+1$
computation, and 863323 times of $x/2$ computation, which is the
largest integer being verified in the world. Wei Ren \cite{weijm}
also pointed out a new approach for the possible proof of Collatz
conjecture. Wei Ren \cite{weiispa} proposed to use a tree-based
graph to reveal two key inner properties in reduced Collatz
dynamics: one is ratio of the count of $x/2$ over the count of
$3x+1$ (for any reduced Collatz dynamics, the count of $x/2$ over
the count of $3x+1$ is larger than $ln3/ln2$), and the other is
partition (all positive integers are partitioned regularly
corresponding to ongoing dynamics). Wei Ren et al. \cite{weihpcc}
also proposed an automata method for fast computing Collatz
dynamics. All source code and output data by computer programs in
those related papers can be accessed in public repository
\cite{weidata}. Jeffrey P. Dumont et al. present beautiful graphs
for the original Collatz dynamics \cite{jeffrey01,jeffrey03} in
2003. The problem is original dynamics mix the characteristics of
reduced dynamics and thus original dynamics presents certain chaos.
Marc Chamberland reviews the works on Collatz conjecture
\cite{marc06} in 2006, and some aspects in available analysis
results are surveyed. Jeffrey C. Lagarias edits a book on 3x+1
problem and also reviews the problem in 2010 \cite{lagarisa10}.

\section{Preliminaries \label{sec:rcc}}

\begin{notation}

\item (1) $[1]_2=\{x|x \equiv 1 \mod 2, x \in Z^+\};$ $[0]_2=\{x|x
\equiv 0 \mod 2, x \in Z^+\}.$

\item (2) $[i]_m=\{x|x \equiv i \mod m, x \in Z^+, m \geq 2, m \in
Z^+, 0 \leq i \leq m-1, i \in  Z^+ \cup \{0\}\}.$

\item (3) $\min(S=\{...\}):$ The minimal element in a set $S$.

\item (4) $\max(S=\{...\}):$ The maximal element in a set $S$.

\end{notation}

In this paper, we only discuss positive integers (denoted as $Z^+$).
Any odd $x$ will iterate to $3x+1$, which is always even. Collatz
computation afterward is always $x/2$. If combing these two as
$(3x+1)/2$, then the Collatz computation $T(x)$ can be defined as
follows: $T(x)=(3x+1)/2$ if $x$ is odd; Otherwise, $T(x)=x/2$. For
the convenience in presentation, we denote $(3x+1)/2$ as `$I(x)$'
(or just `$I$') and $x/2$ as `$O(x)$' (or just `$O$'). Indeed, `$I$'
is named from ``Increase'' due to $(3x+1)/2>x$, and `$O$' is named
from ``dOwn'' due to $x/2<x$.

$T^{(k+1)}(T^{(k)}(x))$ ($k\in Z^+$) means two successive Collatz
computations, where $T^{(k+1)}=I$ if $T^{(k)}(x) \in [1]_2$, and
$T^{(k+1)}=O$ if $T^{(k)}(x) \in [0]_2$. For simplicity by using
less parentheses, we can rewrite it as $T^{(k)}T^{(k+1)}(x)$.

Iteratively, $T^{(k)}(T^{(k-1)}(...(T^{(1)}(x))))$ $k\geq 2, k \in
Z^+$ can be written as $T^{(1)}...T^{(k-1)}T^{(k)}(x)$, and
$T^{(k)}=I$ if $T^{(1)}...T^{(k-1)}(x) \in [1]_2$ and $T^{(k)}=O$ if
$T^{(1)}...T^{(k-1)}(x)\in [0]_2$.

Stopping time of $n \in Z^+$ is defined as the minimal number of
steps needed to iterate to 1:

$$s(n)=inf\{k: T^{(1)}...T^{(k-1)}T^{(k)}(n)=1\}.$$

$T(x)$ is either $(3x+1)/2$ or $x/2$ (i.e., $T\in\{I,O\}$), the
$s(n)$ is thus the count of $(3x+1)/2$ computation plus the count of
$x/2$ computation. 

The Collatz computation sequence (i.e., original dynamics) of $n \in
Z^+$ is the sequence of Collatz computations that occurs from
starting integer to 1:

$$d(n)=T^{(1)}...T^{(L-1)}T^{(L)},$$ where
$T^{(1)}...T^{(L-1)}T^{(L)}(n)=1,$ $L=s(n),$ $T^{(i)} \in \{I,O\},
i=1,...,L.$

For example, Collatz computation sequence from starting integer 3 to
1 is $IIOOO$, because $3 \rightarrow 5 \rightarrow 8 \rightarrow 4
\rightarrow 2 \rightarrow 1.$ Thus, $s(3)=5$. $d(3)=IIOOO.$

\begin{definition} $IsMatched:  x \times t  \rightarrow b$. It takes
as input $x \in Z^+$ and $t \in \{I,O\}$, and outputs $b \in
\{True,False\}.$ If $x \in [1]_2$ and $t=I$, or if $x \in [0]_2$ and
$t=O$, then output $b=True$; Otherwise, output $b=False.$
\end{definition}

Simply speaking, this function checks whether the forthcoming
Collatz transformation (i.e., $t \in \{I,O\}$) matches with the
current transformed integer $x$.

\begin{definition} $Substr: s \times i \times j \rightarrow s'.$ It
takes as input $s,i,j$, where $s \in \{I,O\}^{|s|}$, $1 \leq i \leq
|s|,$ $1\leq j \leq |s|-(i-1)$, and outputs $s'$ where
$s=s_a\|s'\|s_b$, $|s_a|=i-1,$ $|s'|=j$, $|s_b|=|s|-|s_a|-|s'|$ and
``$|\cdot|$'' returns length. \label{df:substr} \end{definition}

\begin{remark}

\item (1) For example, $Substr(IIOO,1,4)=IIOO,
Substr(IIOO,1,3)=IIO.$

\item (2) Especially, $Substr(s,1,|s|)=s.$ $Substr(s,|s|,1)$ returns
the last transformation in $s$. $Substr(s,1,1)$ returns the first
transformation in $s$. $Substr(s,j,1)$ returns the $j$-th
transformation in $s.$

\item (3) In other words, $s'$ is a selected segment in $s$ that
starts from the location $i$ and has the length of $j$.

\item (4) Simply speaking, this function can obtain the Collatz
transforms from $i$ to $i+j-1$ from a given inputting transform
sequence (e.g., reduced dynamics) in terms of $s \in \{I,O\}^{|s|}$.

\item (5) Note that, $Substr(\cdot)$ itself is a function. In other
words, it can be looked as $Substr(\cdot)(\cdot)$. E.g.,
$Substr(IIOO,1,1)(3)=I(3)=(3*3+1)/2=5,$\\
$Substr(IIOO,1,2)(3)=II(3)=I(I(3))=I(5)=(3*5+1)/2=8,$\\
$Substr(IIOO,1,3)(3)=IIO(3)=O(II(3))=O(8)=8/2=4,$\\
$Substr(IIOO,1,4)(3)=IIOO(3)=O(IIO(3))=O(4)=4/2=2<3.$

\item (6) It is worth to stress that, although in above definition
$j \geq 1$, it can be extended to $j \geq 0$ by assuming
$Substr(\cdot,\cdot,0)(x)=x$.

\end{remark}

\begin{definition} Collatz Conjecture.\\ $\forall x \in
\mathbb{Z}^+,$ $\exists L \in \mathbb{Z}^+$, such that
$T^{(1)}...T^{(L-1)}T^{(L)}(x)=1$ where $T^{(j)} \in \{I,O\}$,
$j=1,2,...,L.$ Recall that, \\
$IsMatch(Substr(T^{(1)}...T^{(L-1)}T^{(L)},1,i)(x),$
$Substr(T^{(1)}...T^{(L-1)}T^{(L)},i+1,1))=True,$ where
$i=0,1,2,...,L-1$. \end{definition}

Obviously, Collatz conjecture is held when $x=1$. In the following,
we mainly concern $x \geq 2, x \in Z^+.$

\begin{definition} Reduced Collatz Conjecture. $\forall x \in Z^+, x
\geq 2$, $\exists L \in Z^+$, such that
$T^{(1)}...T^{(L-1)}T^{(L)}(x)<x$ and $T^{(1)}...T^{(i-1)}T^{(i)}(x)
\not <x,$ where $i=1,...,L-1$. $T^{(j)} \in \{I,O\}$ for
$j=1,2,...,L$. Recall that, \\
$IsMatch(Substr(T^{(1)}...T^{(L-1)}T^{(L)},1,i)(x),$
$Substr(T^{(1)}...T^{(L-1)}T^{(L)},i+1,1))=True,$ where
$i=0,1,2,...,L-1$. \end{definition}

The reduced Collatz computation sequence (i.e., reduced dynamics) of
$n \in Z^+$ is the sequence of Collatz computations that occurs from
starting integer to 1:

$$d_r(n)=T^{(1)}...T^{(L-1)}T^{(L)},$$ where
$T^{(1)}...T^{(L-1)}T^{(L)}(n)<n,$
$T^{(1)}...T^{(i-1)}T^{(i)}(n)\not <n,$ $i=1,2,...,L-1$, $T^{(j)}
\in \{I,O\}, j=1,...,L.$

For example, reduced Collatz computation sequence from starting
integer 3 is $IIOO$, because $3 \rightarrow 5 \rightarrow 8
\rightarrow 4 \rightarrow 2.$ Thus, $d_r(3)=IIOO.$

\begin{theorem} \label{pro:ccrcc}
Collatz Conjecture is equivalent to Reduced Collatz Conjecture.
\end{theorem}
The proof is straightforward \cite{weijm}.

Obviously, reduced dynamics is more primitive than original
dynamics, because original dynamics consists of reduced dynamics.
Simply speaking, reduced dynamics are building blocks of original
dynamics.

\begin{remark}

\item (1) Obviously, $d_r(x \in [0]_2)=O$.

\item (2) $d_r(3)=IIOO,$ $d_r(5)=IO,$ $d_r(7)=IIIOIOO,$ $d_r(9)=IO,$
$d_r(11)=IIOIO.$ Indeed, we design computer programs that output all
$d_r(x), \forall x \in [1,99999999]$ \cite{weidata}. From the data
we discover the property - period and its relation to the length (or
the total count of $x/2$ computation) in reduced dynamics - will be
proved in the following of this paper.

\item (3) In fact, we proved some results on $d_r(x)$ for specific
$x$, e.g., $d_r(x \in [1]_4)=IO,$ $d_r(x\in [3]_{16})=IIOO$,
$d_r(x\in [11]_{32})=IIOIO$, and so on. (More reduced dynamics that
can be directly given are discussed in \cite{weiispa}.)

\item (4) $IIOO$ can be denoted in short as $I^2O^2.$ $IIIOIOO$ can
be denoted in short as $I^3OIO^2.$ In other words, we denote
$\underbrace{I...I}_n$ as $I^n$, and we denote
$\underbrace{O...O}_n$ as $O^n$ where $n\in Z^+, n \geq 2$. We also
assume $I^1=I$, $O^1=O$. $I^0=O^0=\emptyset$ means no transformation
occurs.

\item (5) In fact, we formally proved that the ratio exists in any
reduced Collatz dynamics. That is, the count of $x/2$ over the count
of $3x+1$ is larger than $log_23$ \cite{weiratio}. The ratio can
also be observed and verified in our proposed tree-based graph
\cite{weiispa}.


\end{remark}

\begin{example} \label{ex:implication}

If there exists $d_r(x)$ ($x \in Z^+, x \geq 2$), then

\item (1) $s(x)<x$, where $s=d_r(x);$

\item (2) $Substr(s,1,i)(x) \not <x,$ where $i=1,2,...,|s|-1;$

\item (3) $IsMatched(Substr(s,1,j-1)(x),Substr(s,j,1))=True,$ where
$j=1,2,...,|s|$.

\end{example}

\begin{remark}

\item (1) $s(x)$ is the first transformed integer that is less than
the starting integer.

\item (2) $Substr(s,1,i)(x)$ ($i=1,2,...,|s|-1$) are all
intermediate transformed integers (or integers in obit).

\item (3) When $j=1$, we have
$Substr(s,1,j-1)(x)=Substr(s,1,0)(x)=x.$ $Substr(s,j,1)$ is the
first transformation.

\item (4) If $Substr(s,1,j-1)(x)$ ($j=2,...,|s|$) is current
transformed integer, then $Substr(s,j,1)$ is the next intermediate
Collatz transformation.

\end{remark}

\begin{proposition} \label{pro:endwitho} Given $x \in Z^+,$ if there
exists $d_r(x)$, then $d_r(x)$ ends by $O$. \end{proposition}

\begin{proof} Straightforward due to $I(x)=(3x+1)/2>x$. \qed
\end{proof}

\section{Period Theorem \label{sec:period}}


In this section, we will formally prove $d_r(x+2^L)=d_r(x),$
$L=|d_r(x)|,$ $d_r(x) \in \{I,O\}^{\geq 1}$ in this section. Note
that, interestingly, $L$ is indeed the count of $x/2$ computations
in reduced dynamics (as each $I$ includes one $x/2$ and each $x/2$
has one $x/2$).

\subsection{Observations}

\begin{definition} Period. $\min(\{P|d_r(x+P) = d_r(x), x, P \in
Z^+\})$ is called the period of $x$. \end{definition}

\begin{remark}

\item (1) $\forall x \in [0]_2$, $d_r(x+2)=d_r(x)=O$, period
$P=2^{|d_r(x)|}=2^{|O|}=2^1=2$;

\item (2) $\forall x \in [1]_4,$ $d_r(x+4)=d_r(x)=IO,$ period
$P=2^{|d_r(x)|}=2^{|IO|}=2^2=4$.

\item (3) We thus concentrate on $x \in [3]_4$ in the following.

\end{remark}

For easily understanding, we point out two concerns in the
forthcoming proof.

\begin{enumerate}

\item Obviously, $d_r(3+16)=IIOO=d_r(3)$. We \emph{observed} that
during the computing of $IIOO(3+16)$, intermediate transformed
integers (i.e., $I(3+16), II(3+16), IIO(3+16)$) are odd or even, if
and only if $I(3), II(3), IIO(3)$ are odd or even. Besides, $16$ is
the minimal integer to satisfy above requirements.

\item Formally speaking, the parity of $s(x+P)$ and $s(x)$ are
always identical, where $s=Substr(d_r(x),1,i),$ and
$i=1,2,...,|d_r(x)|-1.$ Also, the parity of $x+P$ and $x$ are
identical. That is, the parity sequence during computing for the
reduced dynamics of starting integer $x+P$ is identical with that of
starting integer $x$, which results in the occurred Collatz
transformations of both are exactly identical. Besides, $P$ is the
minimal integer to satisfy above requirements. This is one concern.

\item The other concern is to prove $s(x+P)<x+P$ and $s(x)<x$ where
$s=d_r(x)$; Also, $s(x+P) \not < x+P$ and $s(x) \not < x$ where
$s=Substr(d_r(x),1,i),$ $i=1,2,...,|d_r(x)|-1.$

\end{enumerate}

A new notation $I'(\cdot)$ is introduced hereby to reveal the
relations among $I(x+P)$, $I(x)$ and $I'(P)$. Let $I'(x)=(3x)/2.$

\begin{example} \label{ex:i'}

\item (1) $I(3+16)=(3(3+16)+1)/2=(3*3+1)/2+3*16/2=I(3)+I'(16),$
$I(3)=(3*3+1)/2=5, I'(16)=3*16/2=24 \in [0]_2,$ $5>3$, $5+24>3+16$.
Thus, next transformation for 3+16 and 3 are both $I$.

\item (2) $II(3+16)=I(I(3)+I'(16))=(3*(I(3)+I'(16))+1)/2=(3I(3)+1)/2+3I'(16)/2=II(3)+I'I'(16),$

$II(3)=I(5)=(3*5+1)/2=8, I'I'(16)=I'(24)=3*24/2=36 \in [0]_2,$
$8>3,$ $8+36>(3+16).$ Thus, next transformation for 3+16 and 3 are
both $O$.

\item (3) $IIO(3+16) = O(II(3)+I'I'(16))=IIO(3)+I'I'O(16),$

$IIO(3)=8/2=4$, $I'I'O(16)=36/2 = 18 \in [0]_2$, $4>3,$
$4+18=22>(3+16)$. Thus, next transformation for 3+16 and 3 are both
$O$.

\item (4) $IIOO(3+16) = O(IIO(3)+I'I'O(16))=IIOO(3)+I'I'OO(16)$

$IIOO(3)=4/2=2$, $I'I'OO(16)=18/2=9$, $2<3,$ $2+9=11<(3+16)$. Thus,
reduced dynamics for 3+16 and 3 both end. \end{example}

\begin{remark} \label{remarkI'}

In above example we can \emph{observe} that $I'(16), I'I'(16),
I'I'O(16)$ are always remained \emph{even}. Thus, they do not
influence the resulting next Collatz transformation (i.e., ``$I$''
or ``$O$'') during computing for the reduced dynamics of starting
integer 3. Therefore, the parity of $s(x+16)$ and $s(x)$ are always
identical, where $x=3$, $s=Substr(d_r(x),1,i),$ and
$i=1,2,...,|d_r(x)|-1.$

\end{remark}

For better presentation, we thus introduce two functions as follows:

\begin{definition} $IsEven: x \rightarrow b$. It takes as input $x
\in Z^+$, and outputs $b \in \{True,False\}$, where $b = True$ if $x
\in [0]_2$ and $b = False$ if $x \in [1]_2.$ \end{definition}

\begin{definition} $Replace: s \rightarrow s'$. It takes as input $s
\in \{I,O\}^{\geq 1}$, and outputs $s' \in \{I',O\}^{\geq 1}$, where
$Substr(s',i,1) = I'$ if $Substr(s,i,1)=I$, and $Susbstr(s',i,1) =
O$ if $Substr(s,i,1)=O$, for $i=1,2,...,|s|$.

\end{definition}

\begin{remark}

\item (1) Simply speaking, replacing all ``$I$'' in ``$s$'' respectively by
``$I'$'' will result in ``$s'$''.

\item (2) Obviously, $\forall s \in \{I,O\}^{\geq 1}$, $|s'|=|s|$
where $s'=Replace(s).$

\item (3) By using above introduced functions, we can restate the
reason in Remark \ref{remarkI'} as follows: $Substr(s',1,i) \in
[0]_2$, where $s'=Replace(s)$, $i=1,2,...,|s'|-1,$ thus the parity
of $s(x+16)$ and $s(x)$ are always identical, where $x=3$,
$s=Substr(d_r(x),1,i),$ and $i=1,2,...,|d_r(x)|-1.$ \end{remark}

\subsection{The Proof of Period Theorem}

\begin{lemma} \label{lemma:thefirstct}

If $P \in [0]_2, x \in Z^+$, then $IsEven(x+P)=IsEven(x).$

\end{lemma}
\begin{proof}
Straightforward. Due to $P \in [0]_2$, if $x \in [1]_2,$ then $x+P
\in [1]_2;$ If $x\in [0]_2,$ then $x+P\in[0]_2.$ 
Thus, $IsEven(x+P)=IsEven(x).$ \qed
\end{proof}

\begin{remark}
Above lemma states that if $P \in [0]_2$, the first Collatz
transformation of $x+P$ 
is identical with that of $x$.
\end{remark}



\begin{lemma} \label{lemma:1|STR|=1} $s(x+P)=s(x)+s'(P)$, where $s
\in \{I,O\},$ $s'=Replace(s),$ $x\in Z^+, P \in [0]_2$. \end{lemma}

\begin{proof}

$IsEven(x+P)=IsEven(x)$ because $P\in [0]_2$, due to Lemma
\ref{lemma:thefirstct}. Thus, the first Collatz transformation of
$x+P$ and the first Collatz transformation of $x$ are identical.

(1) Suppose $x \in [1]_2$, so $s=I.$ Thus, $s'=Replace(s)=I'.$

$s(x+P)=I(x+P)=3((x+P)+1)/2=(3x+1)/2+3*P/2=I(x)+I'(P)=s(x)+s'(P).$

(2) Suppose $x \in [0]_2$, so $s=O.$ Thus, $s'=Replace(s)=O.$

$s(x+P)=O(x+P)=(x+P)/2=x/2+P/2=O(x)+O(P)=s(x)+s'(P).$

Summarizing (1) and (2), $s(x+P)=s(x)+s'(P)$.
%
%
%
\qed
\end{proof}

\begin{lemma} \label{lemma:1} (\textbf{Separation Lemma}.) Suppose
$x \in Z^+,$ $s \in \{I,O\}^{\geq 2},$ $s'=Replace(s).$ If
$Substr(s',1,j)(P) \in [0]_2$, $j=0,1,2,...,|s|-1$, then

(1) $IsEven(Substr(s,1,j)(x+P)) = IsEven(Substr(s,1,j)(x))$;

(2) $Substr(s,1,j+1)(x+P)=\\Substr(s,1,j+1)(x)+Substr(s',1,j+1)(P).$
\end{lemma}

\begin{proof}

(1) $j=0$.

(1.1) $Substr(s',1,j)(P) \in [0]_2$. \\
$Substr(s',1,j)(P)=Substr(s',1,0)(P)=P$, thus $P\in [0]_2.$ Thus,
$IsEven(x+P)=IsEven(x)$ due to Lemma \ref{lemma:thefirstct}. Thus,
the intermediate next Collatz transformation of $x+P$ and $x$ are
identical.

(1.2) $Substr(s,1,j+1)(x+P) \\ =Substr(s,1,1)(x+P) \;\;\;\;\;
\because j=0\\ =Substr(s,1,1)(x)+Substr(s',1,1)(P) \;\;\;\;\;
\because$ Lemma \ref{lemma:1|STR|=1}\\
$=Substr(s,1,j+1)(x)+Substr(s',1,j+1)(P).$

%
%
%
%
%
%
%
(2) $j=1$.

(2.1) Due to (1), $Substr(s,1,1)(x+P) =
Substr(s,1,1)(x)+Substr(s',1,1)(P)$.

Besides, $Substr(s',1,1)(P) \in [0]_2$. Thus,

$IsEven(Substr(s,1,1)(x+P)) = IsEven(Substr(s,1,1)(x))$. Thus, the
intermediate next Collatz transformation of $x+P$ and $x$ are
identical.

(2.2) There exists two cases as follows:

(2.2.1) If $Substr(s,1,j+1)=Substr(s,1,j)\|I$, then

$Substr(s,1,j+1)(x+P) \\ =(Substr(s,1,j)\|I)(x+P)\\
=I(Substr(s,1,j)(x+P)) \\ =I(Substr(s,1,j)(x)+Substr(s',1,j)(P))
\;\;\;\;\; \because (1.2) \\
=(3(Substr(s,1,j)(x)+Substr(s',1,j)(P))+1)/2\\
=(3*Substr(s,1,j)(x)+1)/2+3*Substr(s',1,j)(P)/2\\
=I(Substr(s,1,j)(x))+I'(Substr(s',1,j)(P))\\
=(Substr(s,1,j)\|I)(x)+(Substr(s',1,j)\|I')(P)\\
=Substr(s,1,j+1)(x)+Substr(s',1,j+1)(P).$

(2.2.2) If $Substr(s,1,j+1)=Substr(s,1,j)\|O$, then

$Substr(s,1,j+1)(x+P) \\ =Substr(s,1,j)\|O)(x+P)\\
=O(Substr(s,1,j)(x+P))\\ =O(Substr(s,1,j)(x)+Substr(s',1,j)(P)),
\;\;\;\;\; \because (1.2) \\
=(Substr(s,1,j)(x)+Substr(s',1,j)(P))/2\\
=O(Substr(s,1,j)(x))+O(Substr(s',1,j)(P))\\
=(Substr(s,1,j)\|O)(x)+(Substr(s',1,j)\|O)(P)\\
=Substr(s,1,j+1)(x)+Substr(s',1,j+1)(P).$

(Note that, here $j+1=2.$ Recall that ``$\|$'' is concatenation.)

(3) Similarly, $j=2$.
%

Due to (2), $Substr(s,1,2)(x+P)=Substr(s,1,2)(x)+Substr(s',1,2)(P)$.

Besides, $Substr(s',1,2)(P) \in [0]_2$. Thus,

$IsEven(Substr(s,1,j)(x+P))=IsEven(Substr(s,1,j)(x)).$ Thus, the
intermediate next Collatz transformation of $x+P$ and $x$ are
identical.

(Note that, here $j=2$).

Again, we can prove the following similar to (2.2).

$Substr(s,1,j+1)(x+P)=Substr(s,1,j+1)(x)+Substr(s',1,j+1)(P). \;\;$

(Note that, here $j+1=3.$)

(4) Similarly, we can prove $j=3,...,|s'|-1,$ respectively and
especially in an order. \qed
\end{proof}

\begin{remark}
\item (1) Obviously, above conclusion can be extended to include $|s|=1$ by
Lemma \ref{lemma:thefirstct} and Lemma \ref{lemma:1|STR|=1}.

\item (2) Separation Lemma states the sufficient condition (i.e.,
$Substr(s',1,j)(P) \in [0]_2,$ $j=0,1,2,...,|s|-1$) for guaranteeing
that all intermediate parities of transformed integers $x+P$ are
exactly identical with those of $x$.

\item (3) Separation Lemma is general, as $s$ could be either
original dynamics or reduced dynamics of certain $x \in Z^+.$ That
is the reason we give a special name to this lemma for emphasizing
its importance.

\item (4) In above proof, we assume the number of transformations will be $s$. Note that, it
does not influence the conclusion, as we can use condition
``$j=0,1,2,...$'' instead of ``$j=0,1,2,...,|s|-1$'' to omit the
assumption on the number of transformations.

\end{remark}

We next explore how to compute $Substr(s',1,j)(x)$.

\begin{definition} Function $CntI(\cdot).$ $CntI: s
\rightarrow n$. It takes as input $s \in \{I,O\}^{\geq1}$, and
outputs $n \in \mathbb{N}$ that is the count of $I$ in $s$.
\end{definition}


\begin{example} $CntI(IIOO)=2,$ $CntI(III)=3.$ Obviously, the function name
stems from ``Count the number of $I$''.\end{example}

\begin{lemma} \label{lemma:2} Suppose $s\in \{I,O\}^{\geq 1},$
$s'=Replace(s),$ $x \in Z^+$, we have
$$Substr(s',1,j)(x)=\frac{3^{CntI(Substr(s,1,j))}}{2^j}*x, \;\;
j=1,2,...,|s|.$$ \end{lemma}

\begin{proof}
(1) $|s|=1.$ Thus, $j=1.$

(1.1) If $s=I,$ then $s'=Replace(s)=I'.$

$Substr(s',1,j)(x)=Substr(I',1,1)(x)=I'(x)=3*x/2=3^1*x/2^1\\
=3^{CntI(I)}*x/2^{|I|} = 3^{CntI(Substr(s,1,j))}/2^j*x.$

(1.2) If $s=O,$ then $s'=Replace(s)=O.$

$Substr(s',1,j)(x)=O(x)=x/2=3^0*x/2^1=3^{CntI(O)}*x/2^{|O|}\\ =
3^{CntI(Substr(s,1,j))}/2^j*x.$

(2) $|s| \geq 2.$

(2.1) $j=1$.

(2.1.1) If $Substr(s,1,1)=I,$ then

$Substr(s',1,j)(x)=I'(x)=3*x/2=3^1*x/2^1=3^{CntI(I)}*x/2^{|I|}
=3^{CntI(Substr(s,1,j))}/2^j*x.$

(2.1.2) If $Substr(s,1,1)=O,$ then

$Substr(s',1,j)(x)=O(x)=x/2=3^0*x/2^1=3^{CntI(O)}*x/2^{|O|}\\
=3^{CntI(Substr(s,1,j))}/2^j*x.$

(2.2) Iteratively, for $j=1,2,...,|s'|-1$ in an order (recall that
$|s'|=|s|$).

(2.2.1) If $Substr(s,2,1)=I,$ then

$Substr(s',1,j+1)(x)\\ =(Substr(s',1,j)\|I')(x)\\
=I'(Substr(s',1,j)(x))\\ =3*Substr(s',1,j)(x)/2\\
=3*\frac{3^{CntI(Substr(s,1,j))}}{2^j}*x/2 \;\;\; \because (2.1) for
\; j=1, (2.2) for \; j=2,...,|s'|-1\\
=\frac{3^{CntI(Substr(s,1,j))+1}}{2^{j+1}}*x \\
=\frac{3^{CntI(Substr(s,1,j+1))}}{2^{j+1}}*x.  \;\;\; \because
Substr(s,1,j+1)=Substr(s,1,j)\|I$

(2.2.2) If $Substr(s,2,1)=O,$ then

$Substr(s',1,j+1)(x)\\ =(Substr(s',1,j)\|O)(x)\\
=O(Substr(s',1,j)(x))\\ =Substr(s',1,j)(x)/2 \\
=\frac{3^{CntI(Substr(s,1,j))}}{2^j}*x/2 \\
=\frac{3^{CntI(Substr(s,1,j))}}{2^{j+1}}*x \\
=\frac{3^{CntI(Substr(s,1,j+1))}}{2^{j+1}}*x.  \;\;\; \because
Substr(s,1,j+1)=Substr(s,1,j)\|O$ \qed \end{proof}

\begin{remark} \label{remark:ss'} \item (1) Recall that
$s'=Substr(s',1,|s|)$ and $s=Substr(s,1,|s|)$. Thus, when $j=|s|$,
then
$s'(x)=Substr(s',1,|s|)(x)=\frac{3^{CntI(Substr(s,1,|s|))}}{2^{|s|}}*x
=\frac{3^{CntI(s)}}{2^{|s|}}*x$.

\item (2) Indeed, $j=|Substr(s,1,j))|,$ thus the lemma can be
restated as\\
$Substr(s',1,j)(x)=\frac{3^{CntI(Substr(s,1,j))}}{2^{|Substr(s,1,j)|}}*x,$
$j=1,2,...,|s|.$ \end{remark}

\begin{lemma} \label{lemma:3} \[\min(\{P|Substr(s',1,j)(P) \in
[0]_2, j=0,1,2,...,|s|-1,\] \[s \in \{I,O\}^{\geq 2}, s'=Replace(s),
P \in Z^+\})=2^{|s|}.\] \end{lemma}

\begin{proof}

(1) $j=0$, $Substr(s',1,j)(P) \in [0]_2 \Leftrightarrow P \in
[0]_2.$

(2) $j=1,2,...,|s|-1$.

By Lemma \ref{lemma:2},
$Substr(s',1,j)(x)=\frac{3^{CntI(Substr(s,1,j))}}{2^j}*x.$ Thus,

$Substr(s',1,j)(P) =\frac{3^{CntI(Substr(s,1,j))}}{2^j}*P \in [0]_2,
j=1,2,...,|s|-1 \\ \Leftrightarrow P/2^j \in [0]_2,
j=1,2,...,|s|-1\\ \Leftrightarrow P/2^{j+1} \in Z^+,
j=1,2,...,|s|-1\\ \Leftrightarrow P/2^{|s|} \in Z^+\\
\Leftrightarrow \min(P)=2^{|s|}.$

By (1) and (2), $\min(P)=2^{|s|}$.\qed
\end{proof}

%
%

\begin{notation}

$Set_{rd}=\{s|x \in Z^+, \exists d_r(x), s=d_r(x) \in \{I,O\}^{\geq
1}\}.$

\end{notation}

Simply speaking, $Set_{rd}$ is a set of all reduced dynamics for
those $x \in Z^+$ if there exists reduced dynamics of $x$. That is,
$\forall x \in Z^+$, if there exists $d_r(x)$, then $d_r(x)=s$ will
be included in $Set_{rd}$, which is a set of existing reduced
dynamics.

\begin{lemma} \label{lemma:<p} $s'(P) < P,$ where $s \in Set_{rd},$
$s'=Replace(s),$ $P \in Z^+$.
\end{lemma}

\begin{proof} By Lemma \ref{lemma:2} and Remark \ref{remark:ss'}
(1),

$s'(P)=3^{CntI(s)}/2^{|s|}*P.$

Due to Corollary \ref{cor:ratio}, $\frac{3^{CntI(s)}}{2^{|s|}}<1$.
Thus, $3^{CntI(s)}/2^{|s|}*P < P$. \qed
\end{proof}

\begin{remark}
Here we use Corollary \ref{cor:ratio} (\emph{Form Corollary}) given
in Appendix, which is formally proved by us in another paper
\cite{weiratio}.
\end{remark}

\begin{lemma} \label{lemma:not<p} $Substr(s',1,j)(P) > P,
j=1,2,...,|s|-1,$ where $s \in Set_{rd},$ $s'=Replace(s),$ $P \in
Z^+$. \end{lemma}

\begin{proof} By Lemma \ref{lemma:2},

$Substr(s',1,j)(P)=\frac{3^{CntI(Substr(s,1,j))}}{2^j}*P.$

Due to Corollary \ref{cor:ratiosegment},
$\frac{3^{CntI(Substr(s,1,j))}}{2^j}>1.$ Thus
$\frac{3^{CntI(Substr(s,1,j))}}{2^j}*P > P$. \qed \end{proof}


\begin{theorem} \label{th:2^d} (\textbf{Period Theorem}.)

If there exists $d_r(x) \in \{I,O\}^{\geq 1}$, then there exists
$d_r(x+2^L)$, and $d_r(x+2^L)=d_r(x)$ where $L=|d_r(x)|.$

\end{theorem}

\begin{proof}

(1) Regarding two special cases:

$d_r(x\in [0]_2)=O$, $d_r(x+2^L)=d_r(x), L=|d_r(x)|=|O|=1.$

$d_r(x\in [1]_4)=IO$, $d_r(x+2^L)=d_r(x), L=|d_r(x)|=|IO|=2.$

(2) Next, w.l.o.g., suppose $d_r(x)=s, |s| \geq 3.$ Let
$P=2^L=2^{|d_r(x)|}=2^{|s|},$ $s'=Replace(s)$.

Thus, $Substr(s',1,j)(P) \in [0]_2, j = 0,1,2,...,|s|-1$, and $P$ is
the minimal integer for this requirement by Lemma \ref{lemma:3}.

(2.1) Regarding the ordered parity sequence of $x+P$ and $x$.

$P \in[0]_2$, thus the first transformation of $x+P$ and $x$ is
identical by Lemma \ref{lemma:thefirstct}.

$Substr(s,1,j+1)(x+P)=Substr(s,1,j+1)(x)+Substr(s',1,j+1)(P)$ where
$s'=Replace(s)$, by Lemma \ref{lemma:1} (i.e., \emph{Separation
Lemma}). The parity for all transformed integers for $x+P$ and $P$
(except for the last one) are exactly identical due to
$Substr(s',1,j)(P) \in [0]_2, j = 0,1,2,...,|s|-1.$

(2.2) Regarding the comparison between transformed integers and
starting integer.

(2.2.1) $s(x)<x$ due to the definition of reduced dynamics of $x$
(recall that $s=d_r(x)$).

(2.2.2) $Substr(s,1,j)(x) \not <x, j=1,2,...,|s|-1,$ due to the
definition of reduced dynamics of $x$ (recall that $s=d_r(x)$).

(2.2.3) $s'(P)<P$ by Lemma \ref{lemma:<p}.

(2.2.4) $Substr(s',1,j)(P) \not <P, j=1,2,...,|s'|-1,$ by Lemma
\ref{lemma:not<p}.

(2.2.5) $s(x+P)=s(x)+s'(P)<x+P,$ by Lemma \ref{lemma:1}, (2.2.1) and
(2.2.3).

(2.2.6) $Substr(s,1,j)(x+P)= Substr(s,1,j)(x)+ Substr(s',1,j)(P)
\not <x+P, j=1,2,...,|s|-1,$ due to (2.2.2) and (2.2.4).

Thus, only the last transformed integer $s(x+P)$ is less than the
starting integer $x+P$, and the other transformed integers are not
less than the starting integer $x+P$. Thus, $d_r(x+P)=s$.

Due to (2.1) and (2.2), $d_r(x+P)=d_r(x), P=2^L, L=|s| \geq 3.$

Summarizing (1) and (2), $d_r(x+2^L)$ exists and $d_r(x+2^L)=d_r(x)$
where $L=|d_r(x)| \in Z^+.$ \qed \end{proof}

\begin{remark}

Interestingly, the period equals 2 to the power of the count of
$x/2$ in reduced dynamics, as the count of $x/2$ equals to the
length of $d_r(x)$. The count of $x/2$ has two folders: one equals
the count of ``$I$'' due to $(3*x+1)/2$; the other equals the count
of ``$O$'' due to $x/2$. Obviously, the count of ``$I$'' add the
count of ``$O$'' equals the length of $d_r(x)$. \end{remark}

\begin{corollary} \label{cor:1} If reduced dynamics of $x\in Z^+$
exists, denoted as $d_r(x) \in \{I,O\}^L, L\in Z^+$, then reduced
dynamics of $x+k*2^L \; (k \in Z^+)$ exists and is identical with
that of $x$. \end{corollary}

%
%
%

\begin{corollary} \label{th:totalnumberofcode} $\forall s \in
Set_{rd},$ $\|\{x|d_r(x)=s\}\|=\|Z^+\|=\aleph_0$ where
``$\|\cdot\|$'' returns the number of a set. \end{corollary}

\begin{proof} Let $2^{|s|}=P \in Z^+$. $d_r(x+P)=s=d_r(x)$. Thus,
$d_r(x+k*P)=s=d_r(x), k \in Z^+$. A bijective mapping from
$\{x|d_r(x)=s\}$ to $Z^+$ can be created as follows: $x
\leftrightarrow 1, x+k*P \leftrightarrow k+1, k \in Z^+.$ Thus,
$\|\{x|d_r(x)=s\}\|=\|Z^+\|=\aleph_0.$  \qed \end{proof}

\section{Conclusion}

We proved the main theorem called Period Theorem (i.e. Theorem
\ref{th:2^d}) after the preparation of some lemmas, and we also
present some corollaries of the main theorem. The major contribution
of the paper is that we not only discover the period in the reduced
dynamics, but also give what the period is or how to compute the
period. The period for the reduced dynamics of $x$ (i.e., $d_r(x)
\in \{I,O\}^L$) is $2^L$ where $L$ is the length of the reduced
dynamics of $x$ consisting of computations either $(3x+1)/2$
(denoted as $I$ computation) or $x/2$ (denoted as $O$ computation).
That is, $L=|d_r(x)|$. In other words, $L$ is also the counts of
$x/2$ computations, as each $I$ has one $x/2$ computation and each
$O$ has one $x/2$ computation. In short, $d_r(x+2^L)=d_r(x)$ where
$L=|d_r(x)|$ That is, if there exists reduced dynamics of $x$, then
there exists reduced dynamics of $x+2^L$ and these two reduced
dynamics are the same, where $L$ is the length of the reduced
dynamics of $x$. Iteratively, there exists the reduced dynamics of
$x+k*2^L$ where $k\in Z^+$, and all are the same with the reduced
dynamics of $x$.

The power of the Period Theorem is that it can be applied directly
for the proof of Reduced Collatz Conjecture, which is equivalent to
Collatz Conjecture (the equivalence is proved in our paper, see
reference \cite{weijm}). That is, we only need to check partial
integers (that is extremely less than before) for verifying or
proving Collatz conjecture. If one integer $x$ is verified for the
existence of reduced dynamics (namely, $d_r(x)$, which consists of
$I$ and $O$), then all integers related to this integer (i.e.,
$x+k*2^L$, $L=|d_r(x)|$, $k \in Z^+$) are verified. Therefore, what
is left (for future work) is to verify that reduced dynamics of
partial integers exists. If it is true, then Reduced Collatz
Conjecture will be true, and Collatz conjecture will be true too.

\section*{Acknowledgement} The research was financially supported by
the Provincial Key Research and Development Program of Hubei (No.
2020BAB105), Knowledge Innovation Program of Wuhan - Basic Research
(No. 2022010801010197), the Opening Project of Nanchang Innovation
Institute, Peking University (No. NCII2022A02), and National Natural
Science Foundation of China (No. 61972366).

\clearpage
\section*{Appendix}

In our another paper \cite{weiratio}, we proved the following
\emph{Form Corollary} that states the requirements on the count of
``$O$'' and the count of ``$I$'' in any reduced dynamics.

\begin{definition} $CntO: s \rightarrow n$. It
takes as input $s \in \{I,O\}^{\geq 1}$, and outputs $n \in
\mathbb{N}^*$ that is the count of ``$O$'' in $s$.
\end{definition}

E.g., $CntO(IIOO)=2$, $CntO(IO)=1.$

\begin{corollary} \label{th:codetheoremsufficiency} (\textbf{Form
Corollary.}) $\forall s \in \{I,O\}^{\geq1},$ $s \in
Set_{\mathsf{RD}}$, if and only if

(1) $s=O;$ Or,

(2) $CntO(s) = \lceil \log_21.5*CntI(s) \rceil$ and $CntO(s')<
\lceil \log_21.5*CntI(s') \rceil$ \\ where $s'=Substr(s,1,i),$
$i=1,2,...,|s|-1,$ $|s|\geq 2.$ \end{corollary}



Following conclusions are all derived from above \emph{Form
Corollary}.

Following corollary states that $Set_{rd}$ can be constructed by
generating $s \in \{I,O\}^{\geq 1}$ that satisfies aforementioned
requirements instead of by conducting concrete Collatz
transformations for all $x \in Z^+$.

\begin{corollary} \label{cor:code} $Set_{rd} =\{O\} \cup \{s|s \in
\{I,O\}^L, L \in Z^+, L \geq 2, \\ CntO(s) = \lceil
\log_21.5*CntI(s) \rceil, \\ CntO(s') < \lceil \log_21.5*CntI(s')
\rceil, s'=Substr(s,1,i), i=1,2,...,L-1 \}.$ \end{corollary}
\begin{proof} It is straightforward due to Corollary
\ref{th:codetheoremsufficiency}. \qed \end{proof}

Following corollary states the relations between $CntI(s)$ and
$CntO(s)+CntI(s)=|s|$, $s \in Set_{rd}$. Note that, $CntI(s)$ is
indeed equal to the count of $3x+1$ computation, and $|s|$ is indeed
equal to the total count of $x/2$ computation in reduced dynamics.

\begin{corollary} \label{cor:ratio} $s \in \{I,O\}^{\geq 1}, s \in
Set_{rd}$, we have

(1) $|s| \geq \lceil \log_23*CntI(s) \rceil;$ (2)
$3^{CntI(s)}<2^{|s|}.$
\end{corollary}
\begin{proof}

(1) When $s=O$, $CntI(s)=0$, $|s|=1.$ $1> \lceil \log_23*0\rceil=0$
by Corollary \ref{cor:code}.

When $s \neq O,$ $s \in Set_{rd} \Rightarrow |s|=\lceil
\log_21.5*CntI(s) \rceil + CntI(s) =\lceil \log_21.5*CntI(s)+
CntI(s) \rceil = \lceil \log_23*CntI(s)\rceil$ by Corollary
\ref{cor:code}.

In summary, $s \in Set_{rd} \Rightarrow |s| \geq \lceil
\log_23*CntI(s) \rceil,$ and note that ``$>$'' is obtained
\emph{when and only when} $s=O$.

(2) $|s| \geq \lceil \log_23*CntI(s)\rceil \Rightarrow |s| \geq
\log_23^{CntI(s)} \Rightarrow 3^{CntI(s)} \leq 2^{|s|} \Rightarrow
3^{CntI(s)} < 2^{|s|}$. \qed
\end{proof}

\begin{corollary} \label{cor:ratiosegment} $s \in Set_{rd}, s \in
\{I,O\}^{\geq 2}$, we have $3^{CntI(Substr(s,1,j))}>2^j,$
$j=1,2,...,|s|-1.$ \end{corollary} \begin{proof} Let $s' =
Substr(s,1,j), j=1,2,...,|s|-1$. Obviously, $|s'|=j.$

$s \in Set_{rd} \\ \Rightarrow CntO(s') < \lceil \log_21.5*CntI(s')
\rceil \;\;\; \because $ Corollary \ref{th:codetheoremsufficiency}\\
$\Rightarrow CntO(s') < \log_21.5*CntI(s')$ $\;\;\; \because
\log_21.5 \not \in \mathbb{Q}, CntI(s), CntO(s) \in Z^+$\\
$\Rightarrow CntO(s')+CntI(s') < \log_23*CntI(s')\\ \Rightarrow |s'|
< \log_23*CntI(s') \Rightarrow 2^{|s'|}<3^{CntI(s')}\\ \Rightarrow
3^{CntI(Substr(s,1,j))}>2^j.$ \qed \end{proof}


\begin{thebibliography}{00}


\bibitem{UpperboundRecord1}
Tomas Oliveira e Silva, \emph{Maximum excursion and stopping time
record-holders for the 3x+1 problem: computational results},
\newblock{\em Mathematics of Computation}, vol. 68, no. 225, pp. 371-384, 1999.

\bibitem{UpperboundRecord2}
Tomas Oliveira e Silva, \emph{Empirical Verification of the 3x+1 and
Related Conjectures}. In \newblock{\em The Ultimate Challenge: The
3x+1 Problem}, (book edited by Jeffrey C. Lagarias), pp. 189-207,
American Mathematical Society, 2010.



\bibitem{weiuic}
Wei Ren, Simin Li, Ruiyang Xiao and Wei Bi, \emph{Collatz Conjecture
for $2^{100000}-1$ is True - Algorithms for Verifying Extremely
Large Numbers}, \newblock{Proc. of IEEE UIC 2018}, Oct. 2018,
Guangzhou, China, 411-416, 2018

\bibitem{weijm}
Wei Ren, \emph{A New Approach on Proving Collatz Conjecture},
\newblock{Journal of Mathsmatics}, Hindawi, April 2019, ID 6129836,
https://www.hindawi.com/journals/jmath/2019/6129836/.

\bibitem{weiispa}
Wei Ren, \emph{Ratio and Partition are Revealed in Proposed Graph on
Reduced Collatz Dynamics}, \newblock{Proc. of 2019 IEEE Intl Conf on
Parallel \& Distributed Processing with Applications, Big Data \&
Cloud Computing, Sustainable Computing \& Communications, Social
Computing \& Networking (ISPA/BDCloud/SocialCom/SustainCom)}, pp.
474-483, 16-28 Dec. 2019, Ximen, China

\bibitem{weihpcc}
Wei Ren, Ruiyang Xiao, \emph{How to Fast Verify Collatz Conjecture
by Automata}, \newblock{Proc. of IEEE 21st International Conference
on High Performance Computing and Communications; IEEE 17th
International Conference on Smart City; IEEE 5th International
Conference on Data Science and Systems (HPCC/SmartCity/DSS)}, pp.
2720-2729, 10-12 Aug. 2019, Zhangjiajie, China

\bibitem{weiratio}
Wei Ren, \emph{A Reduced Collatz Dynamics Maps to a Residue Class,
and its Count of x/2 over Count of 3*x+1 is larger than ln3/ln2},
\newblock{International Journal of Mathematics and Mathematical Sciences},
Hindawi, Volume2020, Article ID 5946759,
https://doi.org/10.1155/2020/5946759, 2020.

\bibitem{weidata}
Wei Ren, \emph{Reduced Collatz Dynamics Data Reveals Properties for
the Future Proof of Collatz Conjecture},
\newblock{Data}, MDPI, 2019, 4, 89, doi:10.3390/data4020089,
https://www.mdpi.com/2306-5729/4/2/89/pdf.

\bibitem{jeffrey01}
Jeffrey P. Dumont, Clifford A. Reiter, \emph{Visualizing Generalized
3x+1 Function Dynamics},
\newblock{Computers \& Graphics}, vol.25,
no.5, 2001, pp. 883-898

\bibitem{jeffrey03} Jeffrey P. Dumont,
Clifford A. Reiter, \emph{Real Dynamics Of A 3-Power Extension Of
The 3x+1 Function},
\newblock{Dynamics of
Continuous, Discrete and Impulsive Systems: A Mathematical
Analysis}, 10, Dec. 2003, pp. 875-893

\bibitem{marc06} Marc Chamberland, \emph{An Update on the 3x+1 Problem}
\newblock{Butlleti de la Societat Catalana de Matematiques}, 18,
pp.19-45,
https://chamberland.math.grinnell.edu/papers/3x\_survey\_eng.pdf

\bibitem{lagarisa10} Jeffrey C.
Lagarias, \emph{The 3x+1 Problem: an Overview},
\newblock{The Ultimate Challenge: The 3x + 1
Problem}, Edited by Jeffrey C. Lagarias. American Mathematical
Society, Providence, RI, 2010, pp. 3šC29,
https://arxiv.org/pdf/2111.02635.pdf


%
%




\end{thebibliography}
\end{document}